\documentclass[journal,twoside,web,twocolumn]{ieeecolor}
\usepackage{generic}
\usepackage{cite}
\usepackage{amsmath,amssymb,amsfonts}
\usepackage{algorithmic}
\usepackage{graphicx}
\usepackage{textcomp}
\usepackage{mathabx}

\usepackage{amsthm}
\usepackage{optidef}
\usepackage{mathtools}

\theoremstyle{definition}

\newtheorem{assumption}{Assumption}

\theoremstyle{plain}
\newtheorem{theorem}{Theorem}

\newtheorem{lemma}{Lemma}
\newtheorem{probstat}{Problem}

\allowdisplaybreaks

\def\BibTeX{{\rm B\kern-.05em{\sc i\kern-.025em b}\kern-.08em
    T\kern-.1667em\lower.7ex\hbox{E}\kern-.125emX}}
\begin{document}
\title{Discrete-Time Fractional-Order Dynamical Networks Minimum-Energy State Estimation}
\author{Sarthak Chatterjee, Andrea Alessandretti, A. Pedro Aguiar, and S\'ergio Pequito
\thanks{S. Chatterjee is with the Department of Electrical, Computer, and Systems Engineering, Rensselaer Polytechnic Institute, Troy NY, USA (e-mail: chatts3@rpi.edu).}
\thanks{A. Alessandretti is with Magneti Marelli S.p.A., Italy (e-mail: andrea.alessandretti@gmail.com).}
\thanks{A. P. Aguiar is with the Department of Electrical and Computer Engineering, Faculty of Engineering, University of Porto, Porto, Portugal (e-mail: pedro.aguiar@fe.up.pt).}
\thanks{S. Pequito is with the Delft Center for Systems and Control, Delft University of Technology, Delft, The Netherlands (e-mail: \mbox{Sergio.Pequito@tudelft.nl}).}}

\maketitle

\begin{abstract}
Fractional-order dynamical networks are increasingly being used to model and describe processes demonstrating long-term memory or complex interlaced dependencies amongst the spatial and temporal components of a wide variety of dynamical networks. Notable examples include networked control systems or neurophysiological networks which are created using electroencephalographic (EEG) or blood-oxygen-level-dependent (BOLD) data. As a result, the estimation of the states of fractional-order dynamical networks poses an important problem. To this effect, this paper addresses the problem of minimum-energy state estimation for discrete-time \mbox{fractional-order} dynamical networks (DT-FODN), where the state and output equations are affected by an additive noise that is considered to be deterministic, bounded, and unknown. Specifically, we derive the corresponding estimator and show that the resulting estimation error is exponentially input-to-state stable with respect to the disturbances and to a signal that is decreasing with the increase of the accuracy of the adopted approximation model. An illustrative example shows the effectiveness of the proposed method on \mbox{real-world} neurophysiological networks.
\end{abstract}


\section{Introduction}

In a wide variety of dynamical networks, it is often seen that a Markovian dependence of the current state on only the previous state is insufficient to describe the long-term behavior of the considered systems~\cite{moon2008chaotic}. This is due to the fact that \mbox{real-world} networks often demonstrate behaviors in which the current system state is dependent on a combination of several past states or the entire gamut of states seen so far in time. Recent works suggest that discrete-time \mbox{fractional-order} dynamical networks (DT-FODN) evince great success in accurately modeling dynamics that show evidence of nonexponential power-law decay in the dependence of the current state on past states, systems exhibiting long-term memory or fractal properties, or dynamics where there are adaptations in multiple time scales~\cite{lundstrom2008fractional,werner2010fractals,turcott1996fractal,thurner2003scaling,teich1997fractal}. These networks include biological swarms~\cite{west2014networks}, networked control systems~\cite{cao2009distributed,chen2010fractional,ren2011distributed}, and cyber-physical systems~\cite{xuecps} to mention a few. Some of these relationships have also been explored in the context of neurophysiological networks constructed from electroencephalographic (EEG), electrocorticographic (ECoG), or blood-oxygen-level-dependent (BOLD) data~\cite{chatterjee2020fractional,magin2006fractional}.

On the other hand, the problem of state estimation entails the retrieval of the internal state of a given network, often from incomplete or partial measurements of the network's inputs and outputs. Solving this problem is of utmost importance, since, in the majority of real-world networks exchanging measurement information with each other, the network's states are often not directly measurable, and a knowledge of the states is needed to, for example, collectively stabilize the system using state feedback. Given the fundamental nature of the problem, the existence of prior art in the context of state estimation of discrete-time fractional-order systems is no surprise~\cite{sabatier2012observability,sierociuk2006fractional,Safari1,Safari2,miljkovic2017ecg,najar2009discrete,chatterjee2019dealing}.

Nonetheless, in practice, the assumptions in Kalman \mbox{filter-like} formulations are restrictive, as they assume white Gaussian additive process and measurement noises, which implies a uniform prevalence in the power spectrum. Due to this reason, we propose the design of a minimum-energy estimation framework for discrete-time fractional-order networks, where we assume that the state and output equations are affected by an additive disturbance and noise, respectively, that is considered to be deterministic, bounded, and unknown. First proposed by Mortensen~\cite{mortensen1968maximum}, and later refined by Hijab~\cite{hijab1979minimum}, minimum-energy estimators produce an estimate of the system state that is the ``most consistent" with the dynamics and the measurement updates of the system~\cite{fleming1997deterministic,willems2002deterministic,buchstaller2020deterministic,swerling1971modern,bonnabel2014contraction,fagnani1997deterministic,krener2003convergence,krener2018minimum,aguiar2006minimum,hassani2009multiple,pequito2009entropy,alessandretti2011minimum,ha2018cooperative,zamani2013minimum,mceneaney1998robust,haring2020stability}.


In summary, the main contribution of this paper is a minimum-energy estimation procedure to estimate the states of a discrete-time fractional-order dynamical network (\mbox{DT-FODN}). In particular, we prove the exponential \mbox{input-to-state} stability of the estimation error when the aforementioned estimator is used to estimate the states of a \mbox{DT-FODN}. We also provide evidence of the efficacy of our approach via a pedagogical example showing the successful estimation of the states of a neurophysiological network constructed using EEG data.

\emph{Notation:} The symbols $\mathbb{R}, \mathbb{R}^{+}, \mathbb{Z}, \mathbb{N}$, and $\mathbb{N}^{+}$ denote, respectively, the set of reals, positive reals, integers, \mbox{non-negative} integers, and positive integers. Additionally, $\mathbb{R}^n$ and $\mathbb{R}^{n \times m}$ represent the set of column vectors of size $n$ and $n \times m$ matrices with real entries and $I$ denotes an identity matrix of appropriate order. For a given square matrix $M \in \mathbb{R}^{n \times n}$, the notation $M \succeq 0$ (respectively, $M \preceq 0$) indicates that the matrix $M$ is positive semidefinite (respectively, negative semidefinite), i.e., $v^{\mathsf{T}} M v \geq 0$ (respectively, $v^{\mathsf{T}} M v \leq 0$) for any $v \in \mathbb{R}^n$. Further, we use $M^{-\mathsf{T}}$ to denote the inverse of $M^{\mathsf{T}}$. We also write $A \succeq B$ and $A \preceq B$ to mean that the matrix $A-B$ is positive semidefinite and negative semidefinite, respectively. The Euclidean norm is denoted by $\| \cdot \|$.

\section{Problem Formulation}

In this section, we introduce DT-FODN and formulate the minimum-energy state estimation problem for DT-FODN.

\subsection{Discrete-time fractional-order dynamical networks}

Consider a left-bounded sequence $\{ x[k] \}_{k \in \mathbb{Z}}$ over $k$, i.e., with $\limsup\limits_{k \to -\infty} \| x[k] \| < \infty$. Then, for any $\alpha \in \mathbb{R}^{+}$, the \mbox{Gr\"{u}nwald-Letnikov} fractional-order difference is defined as
\begin{equation}
\begin{aligned}
\label{eq:diff_op}
    &\Delta^{\alpha} x[k] \coloneqq \sum_{j=0}^{\infty} c_j^{\alpha} x[k-j],\quad c_j^{\alpha} = (-1)^j \binom{\alpha}{j}, \\ 
    &\quad \binom{\alpha}{j} = \begin{cases} 1 &\mbox{if } j = 0, \\
                    \prod_{i=0}^{j-1} \frac{\alpha-i}{i+1} = \frac{\Gamma (\alpha+1)}{\Gamma (j+1) \Gamma (\alpha-j+1)} & \mbox{if } j > 0, \end{cases}
\end{aligned}
\end{equation}
for all $j \in \mathbb{N}$. The summation in~\eqref{eq:diff_op} is well-defined from the uniform boundedness of the sequence $\{ x[k] \}_{k \in \mathbb{Z}}$ and the fact that $|c^{\alpha}_j| \leq \frac{\alpha^j}{j!}$, which implies that the sequence $\{ c^{\alpha}_j \}_{j \in \mathbb{N}}$ is absolutely summable for any $\alpha \in \mathbb{R}^{+}$~\cite{alessandretti2020finite,sopasakis2017stabilising}.

With the above ingredients, a \emph{discrete-time fractional-order dynamical network with additive disturbance} can be described, respectively, by the state evolution and output equations
\begin{subequations}
\begin{equation}
\label{eq:fos_model}
\sum_{i=1}^{l} A_{i} \Delta^{a_{i}} x[k+1]=\sum_{i=1}^{r} B_{i} \Delta^{b_{i}} u[k]+\sum_{i=1}^{s} G_{i} \Delta^{g_{i}} w[k],
\end{equation}
\begin{equation}
\label{eq:fos_model_op}
z[k] = C'_k x[k] + v'[k],
\end{equation}
\end{subequations}
with the variables $x[k] \in \mathbb{R}^n$, $u[k] \in \mathbb{R}^m$, and $w[k] \in \mathbb{R}^p$ denoting the \emph{state}, \emph{input}, and \emph{disturbance} vectors at time step $k \in \mathbb{N}$, respectively. The scalars $a_i \in \mathbb{R}^{+}$ with $1 \leq i \leq l$, $b_i \in \mathbb{R}^{+}$ with $1 \leq i \leq r$, and $g_i \in \mathbb{R}^{+}$ with $1 \leq i \leq s$ are the \emph{\mbox{fractional-order} coefficients} corresponding, respectively, to the state, the input, and the disturbance. The vectors $z[k], v'[k] \in \mathbb{R}^q$ denote, respectively, the \emph{output} and \emph{measurement disturbance} at time step $k \in \mathbb{N}$. We assume that the (unknown but deterministic) disturbance vectors are bounded as
\begin{equation}
    \| w[k] \| \leq b_w, \| v'[k] \| \leq b_{v'}, \; k \in \mathbb{N},
\end{equation}
for some scalars $b_w, b_{v'} \in \mathbb{R}^{+}$. We also assume that the control input $u[k]$ is known for all time steps $k \in \mathbb{N}$. We denote by $x[0] = x(0)$ the initial condition of the state at time $k = 0$. In the computation of the fractional-order difference, we assume that the system is \emph{causal}, i.e., the state, input, and disturbances are all considered to be zero before the initial time (i.e., $x[k] = 0, u[k] = 0$, and $w[k] = 0$ for all $k < 0$).

With the above ingredients, we seek to solve the following problem in this paper.
\begin{probstat}
\label{prob:main_prob}
Consider the quadratic weighted least-squares objective function
\begin{equation}
\begin{aligned}
&\mathcal{J}\left( x[0], \{ w[i] \}_{i=0}^{N-1}, \{ v'[j] \}_{j=1}^N \right) = \sum_{i=0}^{N-1} w[i]^{\mathsf{T}} Q_i^{-1} w[i] \\ &+ \sum_{j=1}^{N} v'[j]^{\mathsf{T}} R_j^{-1} v'[j] + (x[0] - \hat{x}_0)^{\mathsf{T}} P_0^{-1} (x[0] - \hat{x}_0),
\end{aligned}
\end{equation}
subject to the constraints
\begin{subequations}
\begin{equation}
\label{eq:fos_model_main1}
\sum_{i=1}^{l} A_{i} \Delta^{a_{i}} x[k+1]=\sum_{i=1}^{r} B_{i} \Delta^{b_{i}} u[k]+\sum_{i=1}^{s} G_{i} \Delta^{g_{i}} w[k]
\end{equation}
and
\begin{equation}
\label{eq:fos_model_main2}
z[k] = C'_k x[k] + v'[k],
\end{equation}
\end{subequations}
for some $N \in \mathbb{N}$, with the weighting matrices $Q_i$ $(0 \leq i \leq N-1), R_j$ $(1 \leq j \leq N)$, and $P_0$ chosen to be symmetric and positive definite, and $\hat{x}_0$ chosen to be the \emph{a priori} estimate of the system's initial state. The minimum-energy estimation procedure seeks to solve the following optimization problem
\begin{mini}|l|
  {\scriptstyle \{ x[k] \}_{k=0}^{N}, \{ w[i] \}_{i=0}^{N-1}, \{ v'[j] \}_{j=1}^{N}}{\mathcal{J}\left( x[0], \{ w[i] \}_{i=0}^{N-1}, \{ v'[j] \}_{j=1}^N \right)}{}{}
  \addConstraint{\eqref{eq:fos_model_main1} \: \mathrm{and} \: \eqref{eq:fos_model_main2}},
\label{eq:opt_prob_main1}
\end{mini}
for some $N \in \mathbb{N}$.
\end{probstat}

Additionally, we consider the following mild technical assumption to hold.

\begin{assumption}
\label{assum:sum_inv}
The matrix $\sum_{i=1}^l A_i$ is invertible.
\end{assumption}

\section{Minimum-Energy Estimation for Discrete-Time Fractional-Order Dynamical Networks}

In order to derive the solution to Problem~\ref{prob:main_prob}, we will first start with some alternative formulations of the DT-FODN and relevant definitions that will be used in the sequel. Then, we present the solution in Section~\ref{sec:min_eng} and in Section~\ref{sec:iss} we provide some additional properties of the derived solution, i.e., the exponential input-to-state stability of the estimation error. In Section~\ref{sec:simulations}, we present a practical discussion of the results obtained in the context of DT-FODN. All proofs are relegated to the appendix.

We start by considering a truncation of the last $\mathfrak{v}$ temporal components of~\eqref{eq:fos_model}, which we will refer to as the \mbox{$\mathfrak{v}$-approximation} for the DT-FODN. That being said, we note that using Assumption~\ref{assum:sum_inv}, the DT-FODN model in~\eqref{eq:fos_model} can be equivalently written as
\begin{equation}
    x[k+1]=\sum_{j=1}^{\infty} \check{A}_j x[k-j+1] + \sum_{j=0}^{\infty} \check{B}_j u[k-j] + \sum_{j=0}^{\infty} \check{G}_j w[k-j],
\end{equation}
where $\check{A}_j = -\hat{A}_0^{-1} \hat{A}_j$, $\check{B}_j = \hat{A}_0^{-1} \hat{B}_j$, and $\check{G}_j = \hat{A}_0^{-1} \hat{G}_j$ with $\hat{A}_j = \sum_{i=1}^l A_i c_j^{a_i}$, $\hat{B}_j = \sum_{i=1}^r B_i c_j^{b_i}$, and \mbox{$\hat{G}_j = \sum_{i=1}^s G_i c_j^{g_i}$}. Furthermore, for any positive integer $\mathfrak{v} \in \mathbb{N}^{+}$, the DT-FODN model in~\eqref{eq:fos_model} can be recast as
\begin{subequations}
\label{eq:sys}
\begin{equation}
\label{eq:v_app1}
    \tilde{x}[k+1] = \tilde{A}_{\mathfrak{v}} \tilde{x}[k] + \tilde{B}_{\mathfrak{v}} u[k] + \tilde{G}_{\mathfrak{v}} r[k], \qquad \tilde{x}[0] = \tilde{x}_0,
\end{equation}
\begin{equation}
\label{eq:v_app3}
    y[k+1] = C_{k+1} \tilde{x}[k+1] + v[k+1],
\end{equation}
\end{subequations}
where
\begin{equation}
\label{eq:v_app2}
    r[k] = \sum_{j={\mathfrak{v}}+1}^{\infty} \check{A}_j x[k-j+1] + \sum_{j={\mathfrak{v}}+1}^{\infty} \check{B}_j u[k-j] + \sum_{j=0}^{\infty} \check{G}_j w[k-j],
\end{equation}
with the augmented state vector $\tilde{x}[k] = [ x[k]^{\mathsf{T}}, \ldots, x[k-\mathfrak{v}+1]^{\mathsf{T}}, u[k-1]^{\mathsf{T}},\ldots,u[k-\mathfrak{v}]^{\mathsf{T}} ]^{\mathsf{T}} \in \mathbb{R}^{\mathfrak{v} \times (n+m)}$ and appropriate matrices $\tilde{A}_{\mathfrak{v}}, \tilde{B}_{\mathfrak{v}}$, and $\tilde{G}_{\mathfrak{v}}$, where $\tilde{x}_0 = [x_0^{\mathsf{T}},0,\ldots,0]^{\mathsf{T}}$ denotes the initial condition. The matrices $\tilde{A}_{\mathfrak{v}}$ and $\tilde{B}_{\mathfrak{v}}$ are formed using the terms $\{ \check{A}_j \}_{1 \leq j \leq \mathfrak{v}}$ and $\{ \check{B}_j \}_{1 \leq j \leq \mathfrak{v}}$, while the remaining terms $\{ \check{G}_j \}_{1 \leq j < \infty}$ and the state and input components not included in $\tilde{x}[k]$ are absorbed into the term $\tilde{G}_{\mathfrak{v}} r[k]$. Furthermore, we refer to~\eqref{eq:v_app1} as the \emph{$\mathfrak{v}$-approximation} of the DT-FODN presented in~\eqref{eq:fos_model}.


\subsection{Minimum-energy estimator}
\label{sec:min_eng}

First, let us consider the quadratic weighted least-squares objective function
\begin{equation}
\label{eq:objective}
\begin{aligned}
&\mathcal{J}\left( \tilde{x}[0], \{ r[i] \}_{i=0}^{N-1}, \{ v[j] \}_{j=1}^N \right) = \sum_{i=0}^{N-1} r[i]^{\mathsf{T}} Q_i^{-1} r[i] \\ 
&+ \sum_{j=1}^{N} v[j]^{\mathsf{T}} R_j^{-1} v[j] + (\tilde{x}[0] - \hat{x}_0)^{\mathsf{T}} P_0^{-1} (\tilde{x}[0] - \hat{x}_0),
\end{aligned}
\end{equation}
subject to the constraints
\begin{subequations}
\label{eq:syscon}
\begin{align}
\begin{split}
    \label{eq:syscon1}
    \bar{x}[k+1] &= \tilde{A}_{\mathfrak{v}} \bar{x}[k] + \tilde{B}_{\mathfrak{v}} u[k] + \tilde{G}_{\mathfrak{v}} \bar{r}[k],    
\end{split}\\
\begin{split}
    \label{eq:syscon2}
    y[k+1] &= C_{k+1} \bar{x}[k+1] + \bar{v}[k+1],
\end{split}
\end{align}
\end{subequations}
for some $N \in \mathbb{N}$. The weighting matrices $Q_i$ $(0 \leq i \leq N-1)$ and $R_j$ $(1 \leq j \leq N)$ are chosen to be symmetric and positive definite. The term $\hat{x}_0$ denotes the \emph{a priori} estimate of the (unknown) initial state of the system, with the matrix $P_0$ being symmetric and positive definite.

Subsequently, to construct a minimum-energy estimator for the system~\eqref{eq:sys}, we then consider the weighted \mbox{least-squares} optimization problem
\begin{mini}|l|
  {\scriptstyle \{ \bar{x}[k] \}_{k=0}^{N}, \{ \bar{r}[i] \}_{i=0}^{N-1}, \{ \bar{v}[j] \}_{j=1}^{N}}{\mathcal{J}\left( \tilde{x}[0], \{ r[i] \}_{i=0}^{N-1}, \{ v[j] \}_{j=1}^N \right)}{}{}
  \addConstraint{\eqref{eq:syscon1} \: \text{and} \: \eqref{eq:syscon2}},
\label{eq:opt_prob}
\end{mini}
for some $N \in \mathbb{N}$. The following theorem then certifies the solution of the minimum-energy estimation problem posed in~\eqref{eq:opt_prob}.

\begin{theorem}
\label{thm:soln}
Denote by $\hat{x}[k]$ the state vector that corresponds to the solution of the optimization problem~\eqref{eq:opt_prob}. Then, $\hat{x}[k]$ satisfies the recursion
\begin{equation}
\label{eq:opt_sol}
\begin{aligned}
    &\hat{x}[k+1] = \tilde{A}_{\mathfrak{v}} \hat{x}[k] + \tilde{B}_{\mathfrak{v}} u[k] + K_{k+1} \Big( y[k+1] \\ &- C_{k+1} \Big( \tilde{A}_{\mathfrak{v}} \hat{x}[k] + \tilde{B}_{\mathfrak{v}} u[k] \Big) \Big), \quad 0 \leq k \leq N-1,
\end{aligned}
\end{equation}
with initial conditions specified for $\hat{x}_0$ and $\{ u[j] \}_{j=0}^k$, and with the update equations
\begin{subequations}
\label{eq:filter_updates}
\begin{equation}
    K_{k+1} = M_{k+1} C_{k+1}^{\mathsf{T}} ( C_{k+1} M_{k+1} C_{k+1}^{\mathsf{T}} + R_{k+1} )^{-1},
\end{equation}
\begin{equation}
    M_{k+1} = \tilde{A}_{\mathfrak{v}} P_k \tilde{A}_{\mathfrak{v}}^{\mathsf{T}} + \tilde{G}_{\mathfrak{v}} Q_k \tilde{G}_{\mathfrak{v}}^{\mathsf{T}},
\end{equation}
and
\begin{equation}
\label{eq:P_update}
\begin{aligned}
    &P_{k+1} = (I - K_{k+1} C_{k+1}) M_{k+1} (I - K_{k+1} C_{k+1})^{\mathsf{T}} \\
    &+  K_{k+1} R_{k+1} K_{k+1}^{\mathsf{T}} = (I - K_{k+1} C_{k+1}) M_{k+1},
\end{aligned}
\end{equation}
\end{subequations}
with symmetric and positive definite $P_0$.
\end{theorem}

In Theorem~\ref{thm:soln}, the dynamics of the recursion in~\eqref{eq:opt_sol} (with the initial conditions on $\hat{x}_0$ and the values of $\{ u[j] \}_{j=0}^k$ being known) along with the update equations~\eqref{eq:filter_updates} together solve Problem~\ref{prob:main_prob} completely. It is interesting to note here that the output term $y[k+1]$ presented in~\eqref{eq:syscon2} and~\eqref{eq:opt_sol} is the output of the $\mathfrak{v}$-approximated system~\eqref{eq:sys}, which, in turn, is simply a subset of the outputs $z[k+1]$ obtained from~\eqref{eq:fos_model_op}, truncated $\mathfrak{v}$ time steps in the past, provided $v[k]$ and $C_k$ are formed from the appropriate blocks of $v'[k]$ and $C'_k$ for all $k \in \mathbb{N}$.



In what follows, we show that given the $\mathfrak{v}$-approximation outlined in~\eqref{eq:v_app1}, the evolution of the Lyapunov equation admits a solution over time, by establishing the exponential \mbox{input-to-state} stability of the estimation error.

\subsection{Exponential input-to-state stability of the estimation error}
\label{sec:iss}



In order to prove the exponential input-to-state stability of the minimum-energy estimation error, we need to consider the following mild technical assumptions.


\begin{assumption}
\label{assum:agc_bounds}
There exist constants $\underline{\alpha}, \overline{\alpha}, \beta, \gamma \in \mathbb{R}^{+}$ such that
\begin{equation}
    \underline{\alpha} I \preceq \tilde{A}_{\mathfrak{v}} \tilde{A}_{\mathfrak{v}}^{\mathsf{T}} \preceq \overline{\alpha} I, \quad \tilde{G}_{\mathfrak{v}} \tilde{G}_{\mathfrak{v}}^{\mathsf{T}} \preceq \beta I, \:\: \text{and} \quad C_k^{\mathsf{T}} C_k \preceq \gamma I,
\end{equation}
for all $k \in \mathbb{N}$.
\end{assumption}

First, notice that the \emph{state transition matrix} for the dynamics in~\eqref{eq:v_app1} is given by
\begin{equation}
    \Phi(k,k_0) = \tilde{A}_{\mathfrak{v}}^{(k-k_0)}, \quad \text{with} \quad \Phi(k_0,k_0) = I,
\end{equation}
for all $k \geq k_0 \geq 0$. We also consider the \emph{discrete-time controllability Gramian} associated with the dynamics~\eqref{eq:v_app1} described by
\begin{equation}
\label{eq:contr_gramian}
    W_c(k,k_0) = \sum_{i = k_0}^{k-1} \Phi (k,i+1) \tilde{G}_{\mathfrak{v}} \tilde{G}_{\mathfrak{v}}^{\mathsf{T}} \Phi^{\mathsf{T}} (k,i+1),
\end{equation}
and the \emph{discrete-time observability Gramian} associated with~\eqref{eq:v_app1} to be
\begin{equation}
\label{eq:obsv_gramian}
    W_o(k,k_0) = \sum_{i = k_0+1}^{k} \Phi^{\mathsf{T}} (i,k_0) C_i^{\mathsf{T}} C_i \Phi (i,k_0),
\end{equation}
for $k \geq k_0 \geq 0$. We also make the following assumptions regarding \emph{complete uniform controllability} and \emph{complete uniform observability} of the $\mathfrak{v}$-approximated system in~\eqref{eq:v_app1}.
\begin{assumption}
\label{assum:contr_bounds}
The $\mathfrak{v}$-approximated system~\eqref{eq:v_app1} is completely uniformly controllable, i.e., there exist constants $\delta \in \mathbb{R}^{+}$ and $N_c \in \mathbb{N}^{+}$ such that
\begin{equation}
    W_c (k+N_c,k) \succeq \delta I,
\end{equation}
for all $k \geq 0$.
\end{assumption}

\begin{assumption}
\label{assum:obsv_bounds}
The $\mathfrak{v}$-approximated system~\eqref{eq:v_app1} is completely uniformly observable, i.e., there exist constants \mbox{$\varepsilon \in \mathbb{R}^{+}$} and $N_o \in \mathbb{N}^{+}$ such that
\begin{equation}
    W_o (k+N_o,k) \succeq \varepsilon \Phi^{\mathsf{T}} (k+N_o,k) \Phi(k+N_o,k),
\end{equation}
for all $k \geq 0$.
\end{assumption}
Next, we also present an assumption certifying lower and upper bounds on the weight matrices $Q_k$ and $R_{k+1}$ in~\eqref{eq:objective}.

\begin{assumption}
\label{assum:bounds_matrix}
Without loss of generality, we assume that the weight matrices $Q_k$ and $R_{k+1}$ satisfy
\begin{equation}
    \underline{\vartheta} I \preceq Q_k \preceq \overline{\vartheta} I \quad \text{and} \quad \underline{\rho} I \preceq R_{k+1} \preceq \overline{\rho} I,
\end{equation}
for all $k \geq 0$ and constants $\underline{\vartheta}, \overline{\vartheta}, \underline{\rho}, \overline{\rho} \in \mathbb{R}^+$.
\end{assumption}

\subsubsection{Bounds on the covariance matrix $P_k$}

In this section, we establish lower and upper bounds for the matrix $P_k$, which will be required in Section~\ref{sec:exp_inp_st}, where we use an approach using Lyapunov functions in order to show that the estimation error is exponentially input-to-state stable.

\begin{lemma}
\label{lemma:P1}
Given Assumptions~\ref{assum:agc_bounds} and~\ref{assum:contr_bounds} and the constant $\underline{\pi} \in \mathbb{R}^{+}$, we have that
\begin{equation}
    P_k \succeq \underline{\pi} (N_c) I
\end{equation}
holds for all $k \geq N_c$.
\end{lemma}

\begin{lemma}
\label{lemma:P2}
Given Assumptions~\ref{assum:agc_bounds} and~\ref{assum:obsv_bounds} and the constant $\overline{\pi} \in \mathbb{R}^{+}$, we have that
\begin{equation}
    P_k \preceq \overline{\pi} (N_o) I
\end{equation}
holds for all $k \geq N_o$.
\end{lemma}

\subsubsection{Exponential input-to-state stability of the estimation error}
\label{sec:exp_inp_st}

We start with the minimum-energy estimation error $e[k]$, given by
\begin{equation}
\label{eq:err_defn}
    e[k] = \hat{x}[k] - \tilde{x}[k].
\end{equation}
Next, we certify that the estimation error associated with the \mbox{minimum-energy} estimation process is exponentially input-to-state stable.

\begin{theorem}
\label{thm:main_thm_iss}
Under Assumptions~\ref{assum:agc_bounds},~\ref{assum:contr_bounds}, and~\ref{assum:obsv_bounds}, there exist constants $\sigma, \tau, \chi, \psi \in \mathbb{R}^{+}$ with $\tau < 1$ such that the estimation error $e[k]$ satisfies
\begin{equation}
\label{eq:ISS_theorem}
\begin{aligned}
    \| e[k] \| &\leq \max \Bigg\{ \sigma \tau ^{k-k_0} \| e[k_0] \|,\: \chi \max_{k_o \leq i \leq k-1} \| r[i] \|, \\ &\psi \max_{k_o \leq j \leq k-1} \| v[j+1] \| \Bigg\}
\end{aligned}
\end{equation}
for all $k \geq k_0 \geq \max\{ N_c, N_o \}$.
\end{theorem}

\subsection{Discussion}
\label{sec:discussion}

It is interesting to note that the bound on the estimation error $e[k]$ in~\eqref{eq:ISS_theorem} actually depends on $\| r[i] \|$, where $k_0 \leq i \leq k-1$ for all $i \in \mathbb{N}$. In fact, a distinguishing feature of DT-FODN is the presence of a finite non-zero disturbance term in the \mbox{input-to-state} stability bound of the tracking error when tracking a state other than the origin. This disturbance is dependent on the upper bounds on the non-zero reference state being tracked as well as the input. While the linearity of the \mbox{Gr\"unwald-Letnikov} fractional-order difference operator allows one to mitigate this issue in the case of tracking a non-zero exogenous state by a suitable change of state and input coordinates, this approach is not one we can pursue in this paper, since the state we wish to estimate is unknown. However, it can be shown that as the value of $\mathfrak{v}$ in the $\mathfrak{v}$-approximation increases, the upper bound associated with $\| r[i] \|$ decreases drastically since the $\mathfrak{v}$-approximation gives us progressively better representations of the unapproximated system. This further implies that $\| r[i] \|$ in~\eqref{eq:ISS_theorem} stays bounded, with progressively smaller upper bounds associated with $\| r[i] \|$ (and hence, $\| e[k] \|$) with increasing $\mathfrak{v}$.

Lastly, the estimation error associated with the minimum-energy estimation process in~\eqref{eq:err_defn} is defined in terms of the state of the $\mathfrak{v}$-approximated system $\tilde{x}[k]$. In reality, as detailed above, with larger values of $\mathfrak{v}$, the $\mathfrak{v}$-approximated system approaches the real network dynamics, and thus we obtain an expression for the estimation error with respect to the real system in the limiting case, where the input-to-state stability bound presented in Theorem~\ref{thm:main_thm_iss} holds.

\subsection{Illustrative Example}
\label{sec:simulations}

\begin{figure}[ht]
    \centering
    \includegraphics[width=0.5\textwidth]{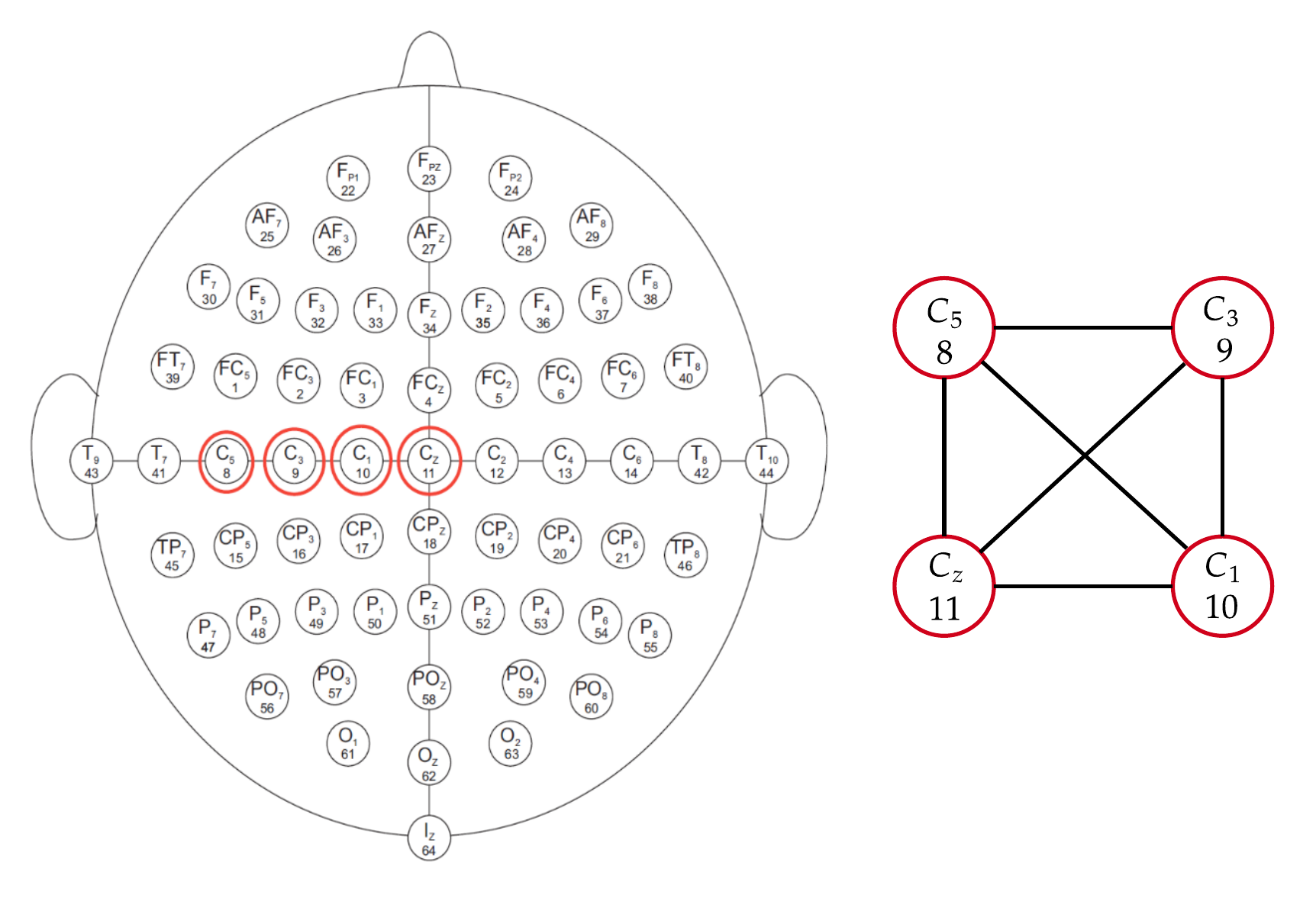}
    \caption{The distribution of the sensors for the measurement of EEG data is shown on the left. The channel labels are shown along with their corresponding numbers and the selected channels over the motor cortex are shown in red. The corresponding network formed by the EEG sensors is shown on the right.}
    \label{fig:eeg_networks}
\end{figure}

In this section, we consider the performance of the minimum-energy estimation paradigm on real-world neurophysiological networks considering EEG data. Specifically, we use $150$ noisy measurements taken from $4$ channels of a $64$-channel EEG signal which records the brain activity of subjects, as shown in Figure~\ref{fig:eeg_networks}. The subjects were asked to perform a variety of motor and imagery tasks, and the specific choice of the $4$ channels was dictated due to them being positioned over the motor cortex of the brain, and, therefore, enabling us to predict motor actions such as the movement of the hands and feet. The data was collected using the BCI$2000$ system with a sampling rate of $160$ Hz~\cite{SchalkBCI,goldberger2000physiobank}. The spatial and temporal parameter components of the DT-FODN assumed to model the original EEG data were identified using the methods described in~\cite{gupta2018dealing}. The matrices $B_i = \begin{bmatrix}1&1&1&1 \end{bmatrix}^{\mathsf{T}}$ for all $i$.

\begin{figure}[ht]
    \centering
    \includegraphics[width=0.5\textwidth]{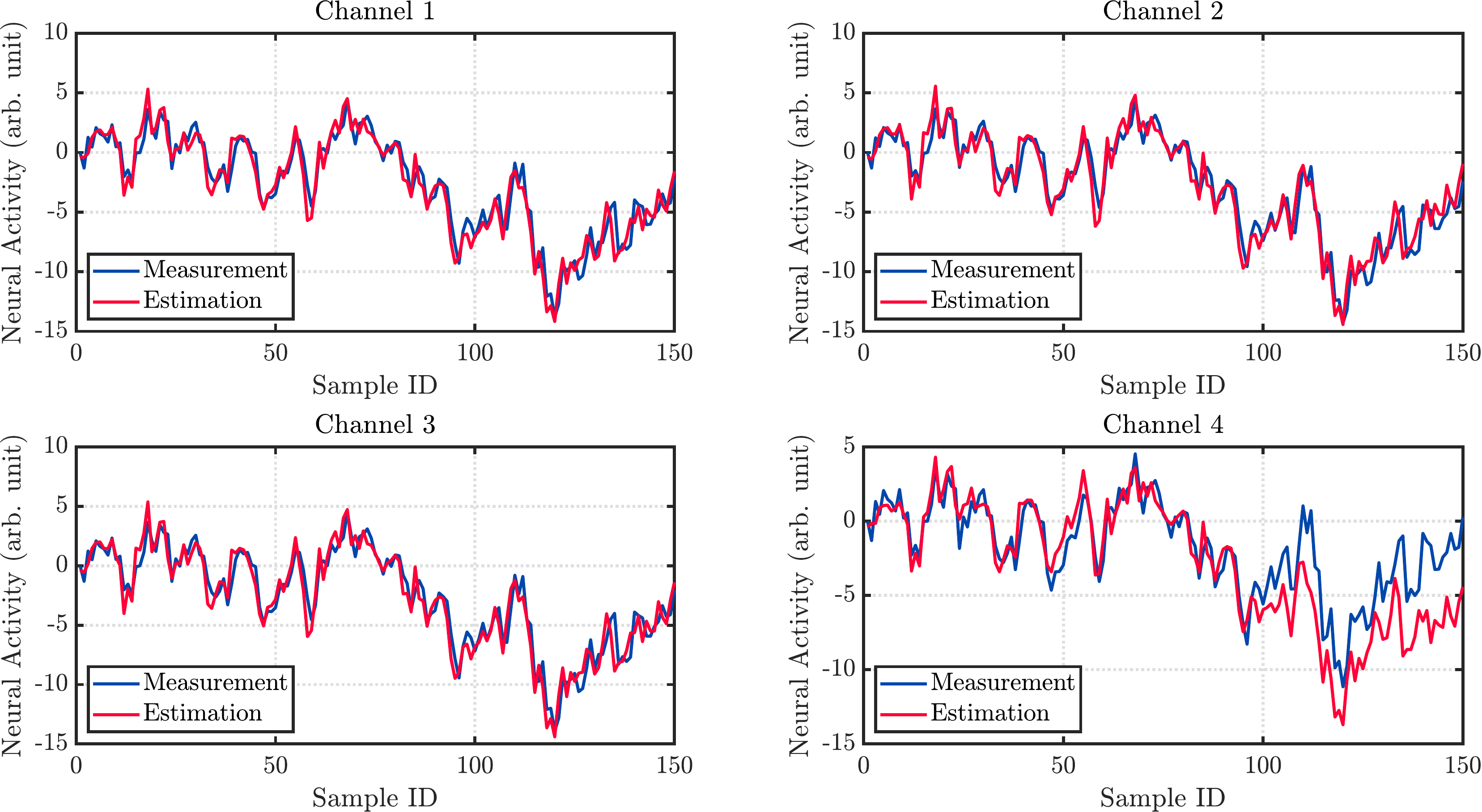}
    \caption{Comparison between the measured output of the $\mathfrak{v}$-augmented system (with $\mathfrak{v}=2$) versus the estimated output of a minimum-energy estimator implemented on the same, in the presence of process and measurement noises for $4$ channels of a $64$-channel EEG signal.}
    \label{fig:output_2}
\end{figure}

\begin{figure}[ht]
    \centering
    \includegraphics[width=0.5\textwidth]{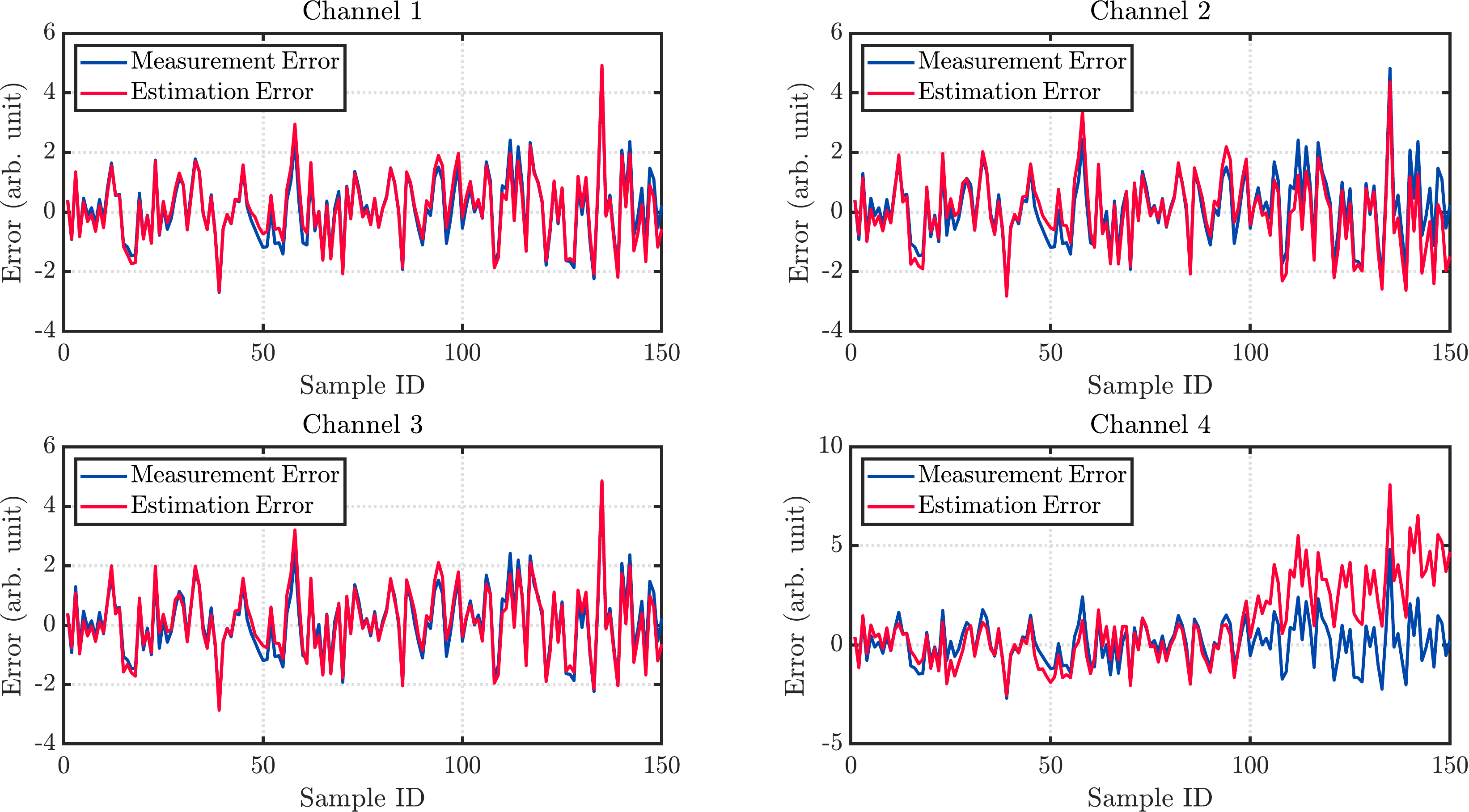}
    \caption{Comparison between the measurement error of the $\mathfrak{v}$-augmented system (with $\mathfrak{v}=2$) versus the estimation error of a minimum-energy estimator implemented on the same, in the presence of process and measurement noises for $4$ channels of a $64$-channel EEG signal.}
    \label{fig:error_2}
\end{figure}

\begin{figure}[ht]
    \centering
    \includegraphics[width=0.5\textwidth]{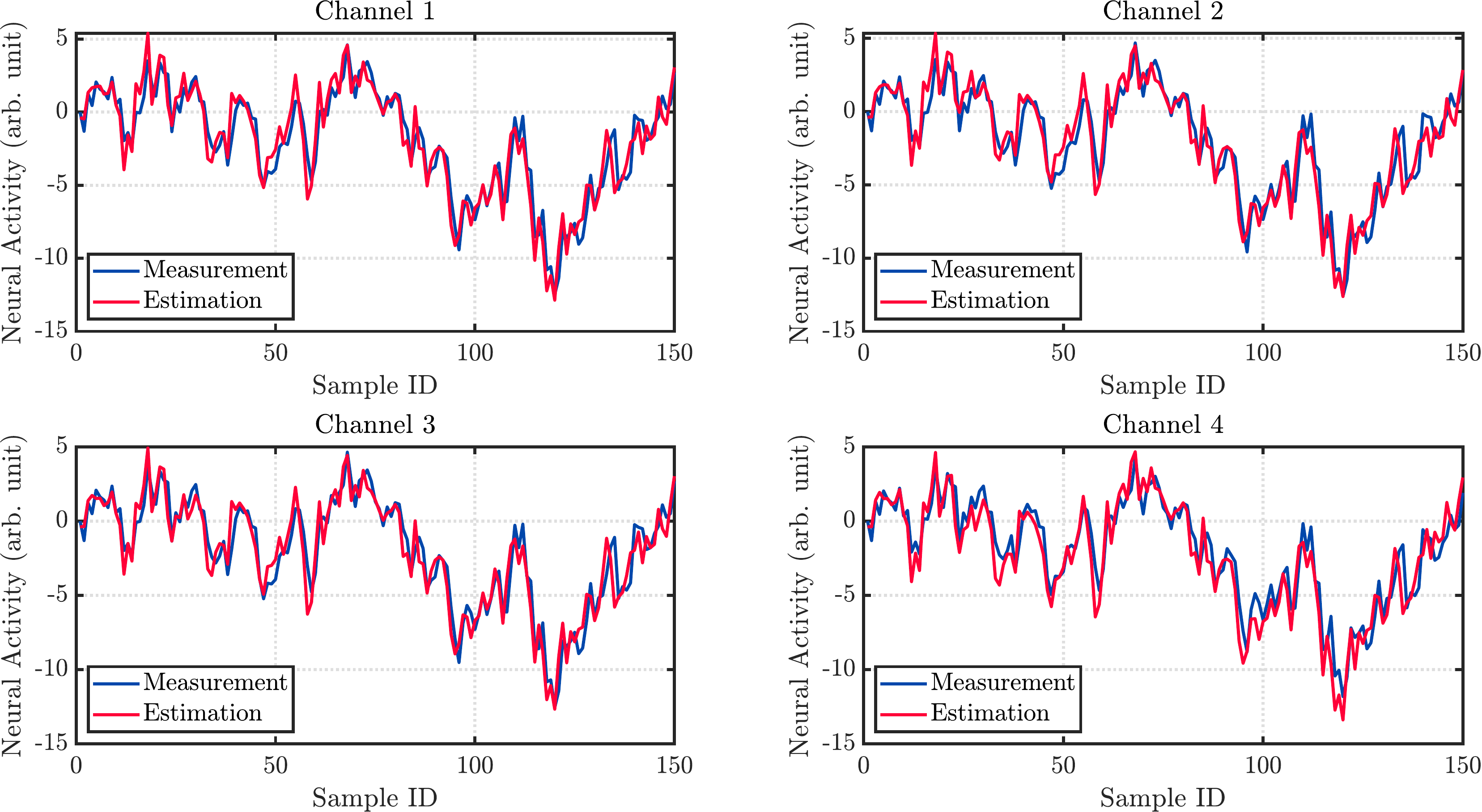}
    \caption{Comparison between the measured output of the $\mathfrak{v}$-augmented system (with $\mathfrak{v}=10$) versus the estimated output of a minimum-energy estimator implemented on the same, in the presence of process and measurement noises for $4$ channels of a $64$-channel EEG signal.}
    \label{fig:output_10}
\end{figure}

\begin{figure}[ht]
    \centering
    \includegraphics[width=0.5\textwidth]{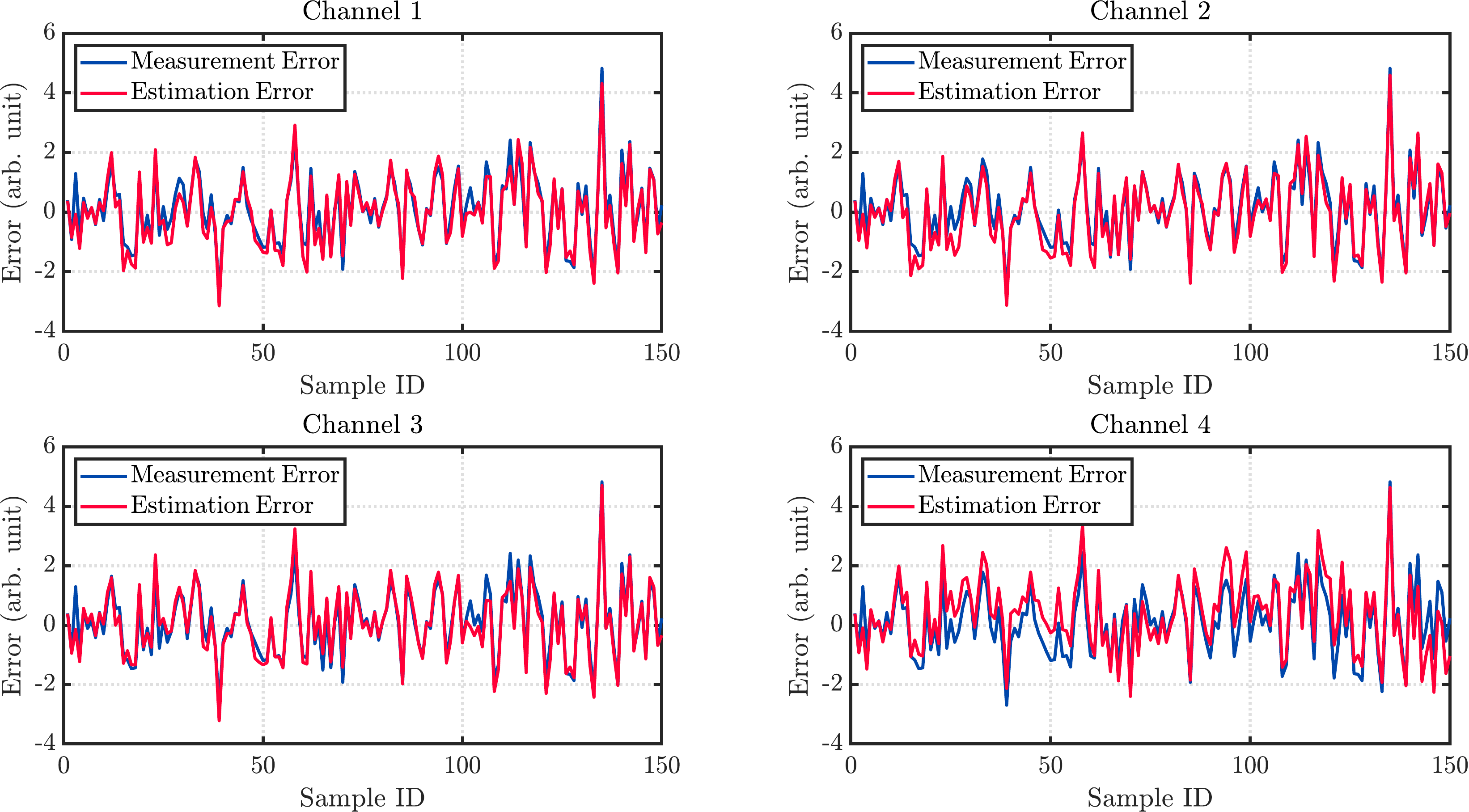}
    \caption{Comparison between the measurement error of the $\mathfrak{v}$-augmented system (with $\mathfrak{v}=10$) versus the estimation error of a minimum-energy estimator implemented on the same, in the presence of process and measurement noises for $4$ channels of a $64$-channel EEG signal.}
    \label{fig:error_10}
\end{figure}

\begin{figure}[ht]
    \centering
    \includegraphics[width=0.5\textwidth]{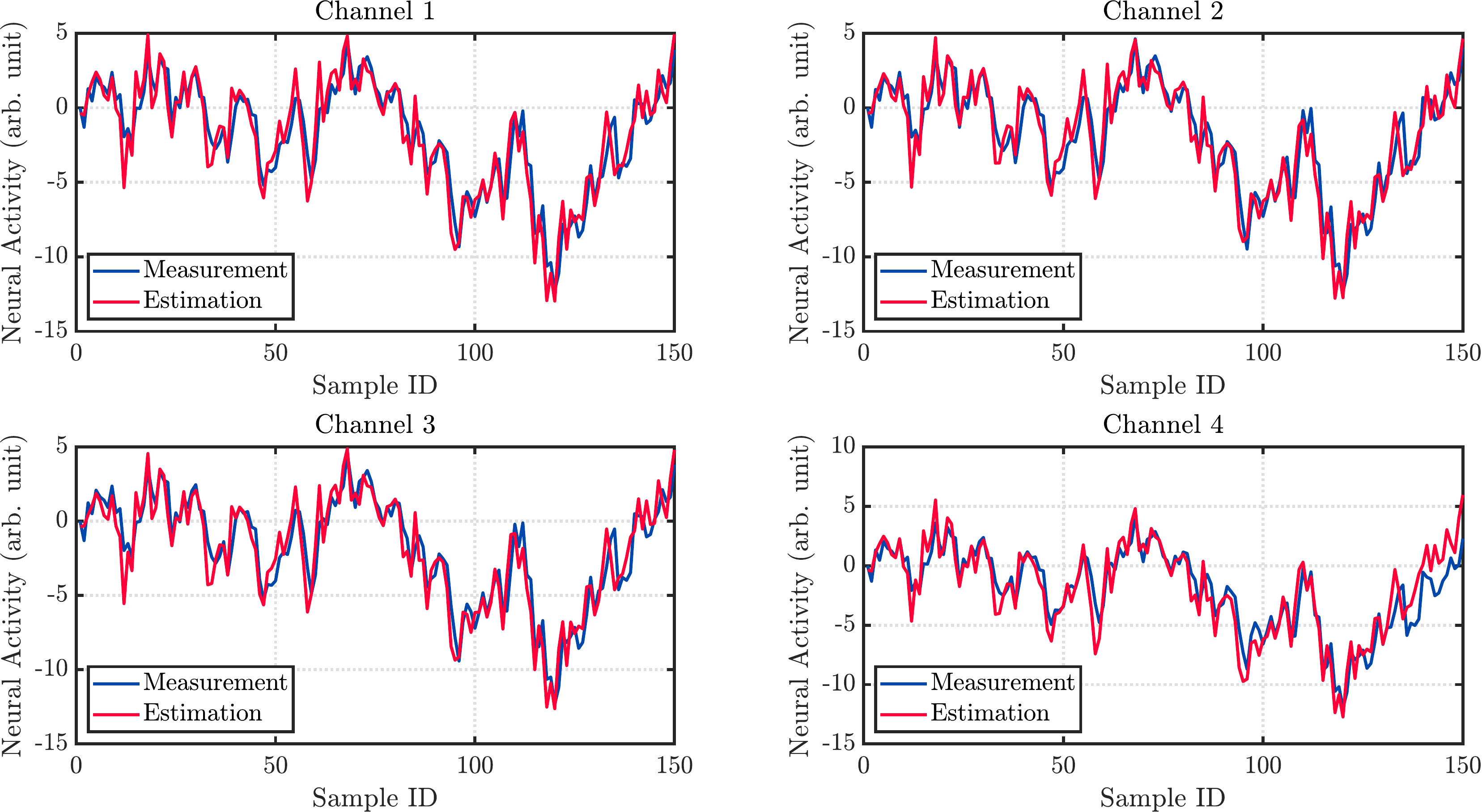}
    \caption{Comparison between the measured output of the $\mathfrak{v}$-augmented system (with $\mathfrak{v}=20$) versus the estimated output of a minimum-energy estimator implemented on the same, in the presence of process and measurement noises for $4$ channels of a $64$-channel EEG signal.}
    \label{fig:output_20}
\end{figure}

\begin{figure}[ht]
    \centering
    \includegraphics[width=0.5\textwidth]{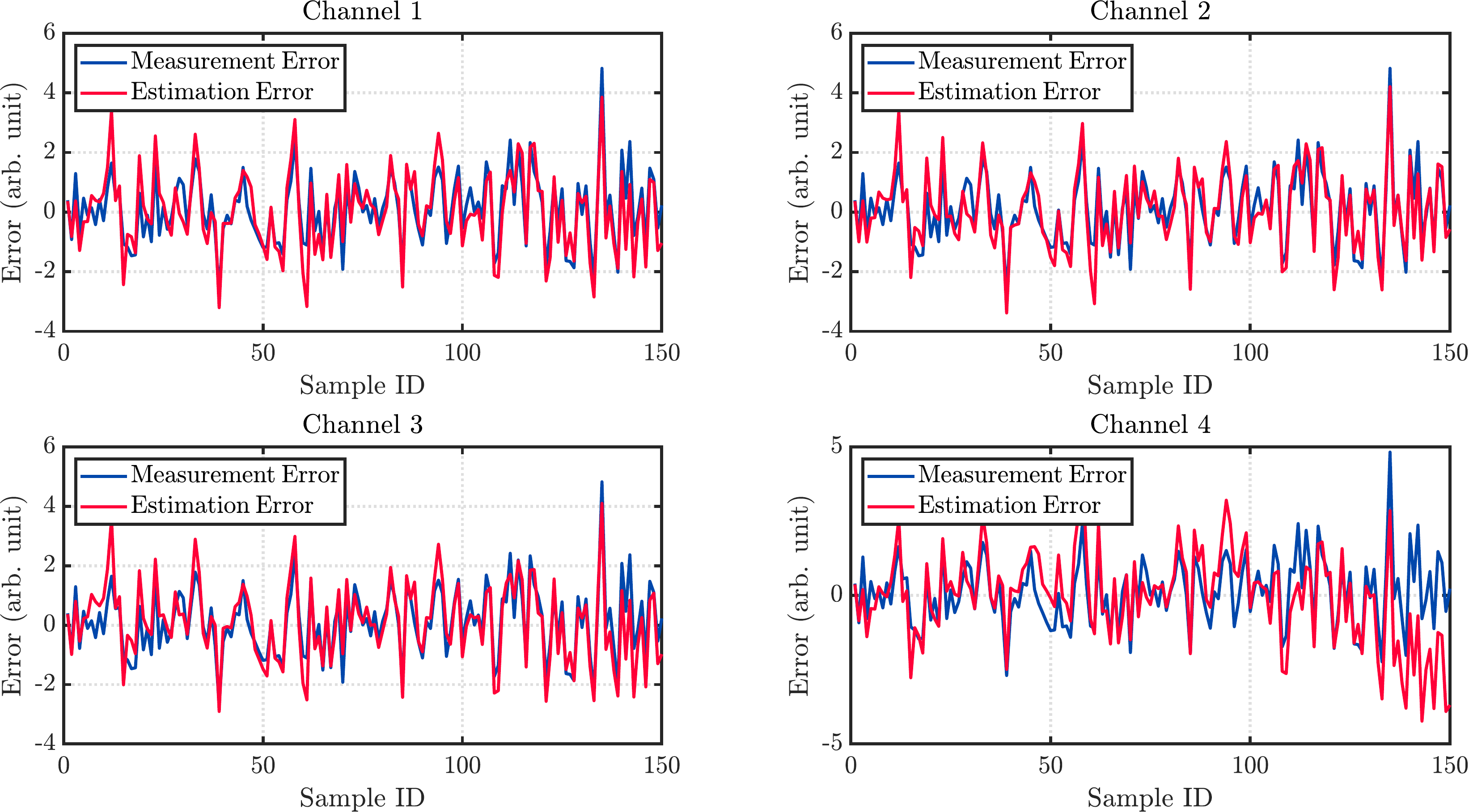}
    \caption{Comparison between the measurement error of the $\mathfrak{v}$-augmented system (with $\mathfrak{v}=20$) versus the estimation error of a minimum-energy estimator implemented on the same, in the presence of process and measurement noises for $4$ channels of a $64$-channel EEG signal.}
    \label{fig:error_20}
\end{figure}

The results of our approach, considering different values of $\mathfrak{v}$, are shown in Figures~\ref{fig:output_2} and~\ref{fig:error_2} (for $\mathfrak{v}=2$), Figures~\ref{fig:output_10} and~\ref{fig:error_10} (for $\mathfrak{v}=10$), and Figures~\ref{fig:output_20} and~\ref{fig:error_20} (for $\mathfrak{v}=20$), which show, respectively (for each value of $\mathfrak{v}$), the comparison between the measured output of the network with noise and the estimated response obtained from the \mbox{minimum-energy} estimator, and also the juxtaposition of the measurement error and the estimation error of the \mbox{minimum-energy} estimation process. We find that the minimum-energy estimator is successfully able to estimate the states in the presence of noise in both the dynamics and the measurement processes.

We also note from the Figures~\ref{fig:output_2} and~\ref{fig:error_2} that when $\mathfrak{v}=2$, we get comparatively larger estimation errors associated with the last $50$ or so samples of Channel $4$, and that this behavior can be mitigated by increasing the value of $\mathfrak{v}$, e.g., by choosing $\mathfrak{v}=10$ or $\mathfrak{v}=20$. This is in line with the discussion in Section~\ref{sec:discussion}, and choosing a larger value of $\mathfrak{v}$ can always, in practice, provide us with better estimation performances, as seen from this example.

\section{Conclusion}

In this paper, we introduced minimum-energy state estimation for discrete-time fractional-order dynamical networks. In particular, the aforementioned minimum-energy estimator is capable of providing an estimate of the unknown states of a \mbox{discrete-time} fractional-order dynamical network while assuming that the associated process and measurement noises are deterministic, bounded, and unknown in nature. We proved that the minimum-energy estimation error is exponentially \mbox{input-to-state} stable and illustrated its performance on \mbox{real-world} neurophysiological EEG networks. Future work will focus on the construction of a resilient and attack-resistant version of the minimum-energy estimator, to take into consideration adversarial attacks or artifacts associated with the measurement process, since the former approach is consistent with the fact that adversarial attacks on sensors often do not follow any particular dynamical or stochastic characterization.

\appendix


\emph{Proof of Theorem~\ref{thm:soln}:} We first consider a single-stage state transition of the system in~\eqref{eq:syscon} and then, sequentially, course through the remaining state transitions. Then, the recursions in~\eqref{eq:opt_sol} and~\eqref{eq:filter_updates} are obtained using the principle of feedback invariance~\cite{hespanha2018linear} and the minimum-energy estimator for discrete-time LTI systems~\cite{buchstaller2020deterministic}, since the $\mathfrak{v}$-approximated \mbox{DT-FODN} in~\eqref{eq:v_app1} fits the latter description.
\qed

\emph{Proof of Lemma~\ref{lemma:P1}:} Suppose $L_{k+1}$ is an arbitrary matrix. We can write
\begin{equation}
\begin{aligned}
    (P_{k+1} + L_{k+1} L_{k+1}^{\mathsf{T}})^{-1} &= \big( ( M_{k+1}^{-1} + C_{k+1}^{\mathsf{T}} R_{k+1}^{-1} C_{k+1} )^{-1} \\
    &+ L_{k+1} L_{k+1}^{\mathsf{T}} \big)^{-1},
\end{aligned}
\end{equation}
where we use the equation
\begin{equation}
\label{eq:mat_inv_P_inv}
    P_{k+1}^{-1} = M_{k+1}^{-1} + C_{k+1}^{\mathsf{T}} R_{k+1}^{-1} C_{k+1},
\end{equation}
which can be obtained from~\eqref{eq:P_update} using the Woodbury identity~\cite[eq.~(157)]{petersen2008matrix}. Notice that the invertibility of $P_k$ and $M_{k+1}$ for any $k \geq 0$ is a consequence of~\eqref{eq:filter_updates}, Assumptions~\ref{assum:agc_bounds} and~\ref{assum:bounds_matrix}, and the fact that $P_0$ is positive definite.

Subsequently, using the bounds in Assumptions~\ref{assum:agc_bounds} and~\ref{assum:bounds_matrix}, and defining $\beta_1 = \frac{\gamma}{\underline{\rho}}$, we have
\begin{align}
\label{eq:big_ineq_P_lower}
    &\left(\mathnormal{P}_{k+1}+\mathnormal{L}_{k+1} \mathnormal{L}_{k+1}^{\mathsf{T}}\right)^{-1} \nonumber \\ &\preceq\left(\left(\mathnormal{M}_{k+1}^{-1}+\beta_{1} \mathnormal{I}\right)^{-1}+\mathnormal{L}_{k+1} \mathnormal{L}_{k+1}^{\mathsf{T}}\right)^{-1} \nonumber \\
    &\stackrel{(\dagger)}{=} \left(\frac{1}{\beta_{1}} \mathnormal{I}-\frac{1}{\beta_{1}^{2}}\left(\mathnormal{M}_{k+1}+\frac{1}{\beta_{1}} \mathnormal{I}\right)^{-1}+\mathnormal{L}_{k+1} \mathnormal{L}_{k+1}^{\mathsf{T}}\right)^{-1} \nonumber \\
    &\stackrel{(\ddagger)}{=}\frac{1}{\beta_{1}^{2}}\left(\frac{1}{\beta_{1}} \mathnormal{I}+\mathnormal{L}_{k+1} \mathnormal{L}_{k+1}^{\mathsf{T}}\right)^{-1} \nonumber \\ & \times \left(\mathnormal{M}_{k+1}+\frac{1}{\beta_{1}} \mathnormal{I}-\frac{1}{\beta_{1}^{2}}\left(\frac{1}{\beta_{1}} \mathnormal{I}+\mathnormal{L}_{k+1} \mathnormal{L}_{k+1}^{\mathsf{T}}\right)^{-1}\right)^{-1} \nonumber \\
    & \times \left(\frac{1}{\beta_{1}} \mathnormal{I}+\mathnormal{L}_{k+1} \mathnormal{L}_{k+1}^{\mathsf{T}}\right)^{-1} +\left(\frac{1}{\beta_{1}} \mathnormal{I}+\mathnormal{L}_{k+1} \mathnormal{L}_{k+1}^{\mathsf{T}}\right)^{-1} \nonumber \\
    &\stackrel{(\diamond)}{=}\frac{1}{\beta_{1}^{2}}\left(\beta_{1} \mathnormal{I}-\beta_{1}^{2} \mathnormal{L}_{k+1}\left(\mathnormal{I}+\beta_{1} \mathnormal{L}_{k+1}^{\mathsf{T}} \mathnormal{L}_{k+1}\right)^{-1} \mathnormal{L}_{k+1}^{\mathsf{T}}\right) \nonumber \\
    & \times \left(\mathnormal{M}_{k+1}+\mathnormal{L}_{k+1}\left(\mathnormal{I}+\beta_{1} \mathnormal{L}_{k+1}^{\mathsf{T}} \mathnormal{L}_{k+1}\right)^{-1} \mathnormal{L}_{k+1}^{\mathsf{T}}\right)^{-1} \nonumber \\
    &\times\left(\beta_{1} \mathnormal{I}-\beta_{1}^{2} \mathnormal{L}_{k+1}\left(\mathnormal{I}+\beta_{1} \mathnormal{L}_{k+1}^{\mathsf{T}} \mathnormal{L}_{k+1}\right)^{-1} \mathnormal{L}_{k+1}^{\mathsf{T}}\right)+\beta_{1} \mathnormal{I} \nonumber \\
    &-\beta_{1}^{2} \mathnormal{L}_{k+1}\left(\mathnormal{I}+\beta_{1} \mathnormal{L}_{k+1}^{\mathsf{T}} \mathnormal{L}_{k+1}\right)^{-1} \mathnormal{L}_{k+1}^{\mathsf{T}}\nonumber\\
    & \stackrel{(\blackdiamond)}{\preceq} 2 \beta_{1}^{2} \mathnormal{L}_{k+1}\left(\mathnormal{I}+\beta_{1} \mathnormal{L}_{k+1}^{\mathsf{T}} \mathnormal{L}_{k+1}\right)^{-1} \mathnormal{L}_{k+1}^{\mathsf{T}} \nonumber \\ & \times \left(\mathnormal{M}_{k+1}+\mathnormal{L}_{k+1}\left(\mathnormal{I}+\beta_{1} \mathnormal{L}_{k+1}^{\mathsf{T}} \mathnormal{L}_{k+1}\right)^{-1} \mathnormal{L}_{k+1}^{\mathsf{T}}\right)^{-1} \nonumber \\
    & \times \mathnormal{L}_{k+1}\left(\mathnormal{I}+\beta_{1} \mathnormal{L}_{k+1}^{\mathsf{T}} \mathnormal{L}_{k+1}\right)^{-1} \mathnormal{L}_{k+1}^{\mathsf{T}} \nonumber \\
    &+2\left(\mathnormal{M}_{k+1}+\mathnormal{L}_{k+1}\left(\mathnormal{I}+\beta_{1} \mathnormal{L}_{k+1}^{\mathsf{T}} \mathnormal{L}_{k+1}\right)^{-1} \mathnormal{L}_{k+1}^{\mathsf{T}}\right)^{-1} \nonumber \\
    &+\beta_{1} \mathnormal{I}-\beta_{1}^{2} \mathnormal{L}_{k+1}\left(\mathnormal{I}+\beta_{1} \mathnormal{L}_{k+1}^{\mathsf{T}} \mathnormal{L}_{k+1}\right)^{-1} \mathnormal{L}_{k+1}^{\mathsf{T}}\nonumber\\
    &\preceq 2(M_{k+1} + \alpha_{1,k+1} L_{k+1} L_{k+1}^{\mathsf{T}})^{-1} + 2 \beta_1 I,
\end{align}
where $\alpha_{1,k+1} = \| I + \beta_1 L_{k+1}^{\mathsf{T}} L_{k+1} \|^{-1}$. The equalities $(\dagger)$, $(\ddagger)$, and $(\diamond)$ in~\eqref{eq:big_ineq_P_lower} are obtained via three successive applications of the Woodbury identity and the inequality $(\blackdiamond)$ in~\eqref{eq:big_ineq_P_lower} is obtained by using the Young-like inequality
\begin{equation}
\label{eq:younglike}
    (f(v) + g(v))^{\mathsf{T}} (f(v)+g(v)) \leq 2 f^{\mathsf{T}}(v) g(v) + 2 g^{\mathsf{T}}(v) f(v),
\end{equation}
with $f(v) = (M_{k+1} + L_{k+1} (I + \beta_1 L_{k+1}^{\mathsf{T}} L_{k+1})^{-1} L_{k+1}^{\mathsf{T}} )^{-\frac{1}{2}} v$ and $g(v) = -\beta_1 (M_{k+1} + L_{k+1} (I + \beta_1 L_{k+1}^{\mathsf{T}} L_{k+1})^{-1} L_{k+1}^{\mathsf{T}} )^{-\frac{1}{2}} L_{k+1} (I + \beta_1 L_{k+1}^{\mathsf{T}} L_{k+1})^{-1} L_{k+1}^{\mathsf{T}} v$.

Plugging in the value of $M_{k+1}$ from the update equations~\eqref{eq:filter_updates}, we have
\begin{equation}
\label{eq:res1}
\begin{aligned}
    &(\mathnormal{P}_{k+1}+\mathnormal{L}_{k+1} \mathnormal{L}_{k+1}^{\mathsf{T}})^{-1} \preceq 2 \beta_1 I + 2 \tilde{A}_{\mathfrak{v}}^{-\mathsf{T}} \\
    \times &\left( P_k + \tilde{A}_{\mathfrak{v}}^{-1} ( \tilde{G}_{\mathfrak{v}} Q_k \tilde{G}_{\mathfrak{v}}^{\mathsf{T}} + \alpha_{1,k+1} L_{k+1} L_{k+1}^{\mathsf{T}} ) \tilde{A}_{\mathfrak{v}}^{-\mathsf{T}} \right)^{-1} \tilde{A}_{\mathfrak{v}}^{-1}.
\end{aligned}
\end{equation}
Now, for any $k \geq 0$, define recursively
\begin{equation}
\label{eq:res2}
    L_j L_j^{\mathsf{T}} = \tilde{A}_{\mathfrak{v}}^{-1} ( \tilde{G}_{\mathfrak{v}} Q_j \tilde{G}_{\mathfrak{v}}^{\mathsf{T}} + \alpha_{1,j+1} L_{j+1} L_{j+1}^{\mathsf{T}} ) \tilde{A}_{\mathfrak{v}}^{-\mathsf{T}}
\end{equation}
for $k \leq j \leq k+N_c-1$, with $L_{k+N_c} L_{k+N_c}^{\mathsf{T}} = 0$. By substituting~\eqref{eq:res2} into~\eqref{eq:res1}, and repeatedly applying the resulting inequality we obtain
\begin{equation}
\label{eq:int1}
\begin{aligned}
&\mathnormal{P}_{k+N_{c}}^{-1} \preceq 2^{N_{c}} \Phi^{-\mathsf{T}}\left(k+N_{c}, k\right)\left(\mathnormal{P}_{k}+\mathnormal{L}_{k} \mathnormal{L}_{k}^{\mathsf{T}}\right)^{-1}\\
&\times \Phi^{-1}\left(k+N_{c}, k\right)+2 \beta_{1} \sum_{i=0}^{N_{c}-1} 2^{i} \\
&\times \Phi^{-\mathsf{T}}\left(k+N_{c}, k+N_{c}-i\right) \Phi^{-1}\left(k+N_{c}, k+N_{c}-i\right).
\end{aligned}
\end{equation}
Using the bounds defined in Assumption~\ref{assum:bounds_matrix},~\eqref{eq:contr_gramian}, and~\eqref{eq:res2}, we can write
\begin{equation}
\label{eq:int2}
    L_k L_k^{\mathsf{T}} \succeq \gamma_1 \Phi^{-1} (k+N_c,k) W_c(k+N_c,k) \Phi^{-\mathsf{T}} (k+N_c,k),
\end{equation}
with $\gamma_1 = \underline{\vartheta} \prod_{j=k}^{k+N_c-1} \alpha_{1,j+1}$. Aggregating the bounds in~\eqref{eq:int1},~\eqref{eq:int2}, and invoking Assumptions~\ref{assum:agc_bounds} and~\ref{assum:contr_bounds}, we have
\begin{equation}
    P_{k+N_c} \succeq \underbrace{\left( \frac{2^{N_c}}{\gamma_1 \delta} + 2 \beta_1 \sum_{i=0}^{N_c-1} \left( \frac{2}{\underline{\alpha}} \right)^i \right)^{-1}}_{\underline{\pi}(N_c)} I.
\end{equation}
\qed


\emph{Proof of Lemma~\ref{lemma:P2}:} Suppose $Y_{k+1}$ is an arbitrary matrix. We can write
\begin{equation}
\begin{aligned}
    (P_{k+1}^{-1} + Y_{k+1}^{\mathsf{T}} Y_{k+1})^{-1} &= \Big( ( \tilde{A}_{\mathfrak{v}} P_k \tilde{A}_{\mathfrak{v}}^{\mathsf{T}} + \tilde{G}_{\mathfrak{v}} Q_k \tilde{G}_{\mathfrak{v}}^{\mathsf{T}} )^{-1} \\
    &+ Z_{k+1}^{\mathsf{T}} Z_{k+1} \Big)^{-1},
\end{aligned}
\end{equation}
where the matrix $Z_{k+1}$ is defined as 
\begin{equation}
\label{eq:Z_bound}
    Z_{k+1} = C_{k+1}^{\mathsf{T}} R_{k+1}^{-1} C_{k+1} + Y^{\mathsf{T}}_{k+1} Y_{k+1}.
\end{equation}
Using the bounds in Assumptions~\ref{assum:agc_bounds} and~\ref{assum:bounds_matrix}, and defining $\beta_2 = \beta \overline{\vartheta}$, we have
\begin{align}
\label{eq:res21}
    &\left(\mathnormal{P}_{k+1}^{-1}+\mathnormal{Y}_{k+1}^{\mathsf{T}} \mathnormal{Y}_{k+1}\right)^{-1} \nonumber \\
    &\preceq\left(\left(\tilde{A}_{\mathfrak{v}} \mathnormal{P}_{k} \tilde{A}_{\mathfrak{v}}^{\mathsf{T}}+\beta_{2} \mathnormal{I}\right)^{-1}+\mathnormal{Z}_{k+1}^{\mathsf{T}} \mathnormal{Z}_{k+1}\right)^{-1} \nonumber \\
    &\stackrel{(\bigtriangleup)}{=} \Bigg(\frac{1}{\beta_{2}} \mathnormal{I}-\frac{1}{\beta_{2}^{2}} \tilde{A}_{\mathfrak{v}}\left(\mathnormal{P}_{k}^{-1}+\frac{1}{\beta_{2}} \tilde{A}_{\mathfrak{v}}^{\mathsf{T}} \tilde{A}_{\mathfrak{v}}\right)^{-1} \tilde{A}_{\mathfrak{v}}^{\mathsf{T}} \nonumber \\
    &+\mathnormal{Z}_{k+1}^{\mathsf{T}} \mathnormal{Z}_{k+1}\Bigg)^{-1} \nonumber \\
    & \stackrel{(\bigtriangledown)}{=} \frac{1}{\beta_{2}^{2}}\left(\frac{1}{\beta_{2}} \mathnormal{I}+\mathnormal{Z}_{k+1}^{\mathsf{T}} \mathnormal{Z}_{k+1}\right)^{-1} \tilde{A}_{\mathfrak{v}} \nonumber \\ 
    &\times \left(\mathnormal{P}_{k}^{-1}+\frac{1}{\beta_{2}} \tilde{A}_{\mathfrak{v}}^{\mathsf{T}} \tilde{A}_{\mathfrak{v}}-\frac{1}{\beta_{2}^{2}} \tilde{A}_{\mathfrak{v}}^{\mathsf{T}}\left(\frac{1}{\beta_{2}} \mathnormal{I}+\mathnormal{Z}_{k+1}^{\mathsf{T}} \mathnormal{Z}_{k+1}\right)^{-1} \tilde{A}_{\mathfrak{v}}\right)^{-1} \nonumber \\
    & \times \tilde{A}_{\mathfrak{v}}^{\mathsf{T}} \left(\frac{1}{\beta_{2}} \mathnormal{I}+\mathnormal{Z}_{k+1}^{\mathsf{T}} \mathnormal{Z}_{k+1}\right)^{-1}+\left(\frac{1}{\beta_{2}} \mathnormal{I}+\mathnormal{Z}_{k+1}^{\mathsf{T}} \mathnormal{Z}_{k+1}\right)^{-1} \nonumber \\
    & \stackrel{(\centerdot)}{=} \frac{1}{\beta_{2}^{2}}\left(\beta_{2} \mathnormal{I}-\beta_{2}^{2} \mathnormal{Z}_{k+1}^{\mathsf{T}}\left(\mathnormal{I}+\beta_{2} \mathnormal{Z}_{k+1} \mathnormal{Z}_{k+1}^{\mathsf{T}}\right)^{-1} \mathnormal{Z}_{k+1}\right) \tilde{A}_{\mathfrak{v}} \nonumber \\
    &\times \left(\mathnormal{P}_{k}^{-1}+\tilde{A}_{\mathfrak{v}}^{\mathsf{T}} \mathnormal{Z}_{k+1}^{\mathsf{T}}\left(\mathnormal{I}+\beta_{2} \mathnormal{Z}_{k+1} \mathnormal{Z}_{k+1}^{\mathsf{T}}\right)^{-1} \mathnormal{Z}_{k+1} \tilde{A}_{\mathfrak{v}}\right)^{-1} \tilde{A}_{\mathfrak{v}}^{\mathsf{T}} \nonumber \\
    &\times\left(\beta_{2} \mathnormal{I}-\beta_{2}^{2} \mathnormal{Z}_{k+1}^{\mathsf{T}}\left(\mathnormal{I}+\beta_{2} \mathnormal{Z}_{k+1} \mathnormal{Z}_{k+1}^{\mathsf{T}}\right)^{-1} \mathnormal{Z}_{k+1}\right) \nonumber \\
    &+\beta_{2} \mathnormal{I}-\beta_{2}^{2} \mathnormal{Z}_{k+1}^{\mathsf{T}}\left(\mathnormal{I}+\beta_{2} \mathnormal{Z}_{k+1} \mathnormal{Z}_{k+1}^{\mathsf{T}}\right)^{-1} \mathnormal{Z}_{k+1} \nonumber\\
    & \stackrel{(\centerdot \centerdot)}{\preceq} 2 \beta_{2}^{2} \mathnormal{Z}_{k+1}^{\mathsf{T}}\left(\mathnormal{I}+\beta_{2} \mathnormal{Z}_{k+1} \mathnormal{Z}_{k+1}^{\mathsf{T}}\right)^{-1} \mathnormal{Z}_{k+1} \tilde{A}_{\mathfrak{v}} \nonumber \\
    & \times \left(\mathnormal{P}_{k}^{-1}+\tilde{A}_{\mathfrak{v}}^{\mathsf{T}} \mathnormal{Z}_{k+1}^{\mathsf{T}}\left(\mathnormal{I}+\beta_{2} \mathnormal{Z}_{k+1} \mathnormal{Z}_{k+1}^{\mathsf{T}}\right)^{-1} \mathnormal{Z}_{k+1} \tilde{A}_{\mathfrak{v}}\right)^{-1} \nonumber \\
    &\times \tilde{A}_{\mathfrak{v}}^{\mathsf{T}} \mathnormal{Z}_{k+1}^{\mathsf{T}} \left(\mathnormal{I}+\beta_{2} \mathnormal{Z}_{k+1} \mathnormal{Z}_{k+1}^{\mathsf{T}}\right)^{-1} \mathnormal{Z}_{k+1} \nonumber \\
    &+2 \tilde{A}_{\mathfrak{v}}\left(\mathnormal{P}_{k}^{-1}+\tilde{A}_{\mathfrak{v}}^{\mathsf{T}} \mathnormal{Z}_{k+1}^{\mathsf{T}}\left(\mathnormal{I}+\beta_{2} \mathnormal{Z}_{k+1} \mathnormal{Z}_{k+1}^{\mathsf{T}}\right)^{-1} \mathnormal{Z}_{k+1} \tilde{A}_{\mathfrak{v}}\right)^{-1} \tilde{A}_{\mathfrak{v}}^{\mathsf{T}}\nonumber\\
    &+\beta_{2} \mathnormal{I} -\beta_{2}^{2} \mathnormal{Z}_{k+1}^{\mathsf{T}}\left(\mathnormal{I}+\beta_{2} \mathnormal{Z}_{k+1} \mathnormal{Z}_{k+1}^{\mathsf{T}}\right)^{-1} \mathnormal{Z}_{k+1} \nonumber\\
    &\preceq 2 \tilde{A}_{\mathfrak{v}} (P_k^{-1} + \alpha_{2,k+1} \tilde{A}_{\mathfrak{v}}^{\mathsf{T}} Z_{k+1}^{\mathsf{T}} Z_{k+1} \tilde{A}_{\mathfrak{v}})^{-1} + 2 \beta_2 I
\end{align}
where $\alpha_{2,k+1} = \| I + \beta_2 Z_{k+1} Z_{k+1}^{\mathsf{T}} \|^{-1}$. The equalities $(\bigtriangleup)$, $(\bigtriangledown)$, and $(\centerdot)$ in~\eqref{eq:res21} are obtained via three successive applications of the Woodbury identity and the inequality $(\centerdot \centerdot)$ in~\eqref{eq:res21} is obtained by using the Young-like inequality~\eqref{eq:younglike}.

Now, for any $k \geq 0$, define
\begin{equation}
\label{eq:res22}
    Y_j^{\mathsf{T}} Y_j = \alpha_{2,j+1} \tilde{A}_{\mathfrak{v}}^{\mathsf{T}} Z_{j+1}^{\mathsf{T}} Z_{j+1} \tilde{A}_{\mathfrak{v}},
\end{equation}
where $k \leq j \leq k+N_o-1$, with $Y_{k+N_o}^{\mathsf{T}} Y_{k+N_o} = 0$. By repeatedly applying~\eqref{eq:res21} and~\eqref{eq:res22}, we obtain
\begin{equation}
\label{eq:int1prime}
\begin{aligned}
&\mathnormal{P}_{k+N_{o}} \preceq 2^{N_{o}} \Phi\left(k+N_{o}, k\right)\left(\mathnormal{P}_{k}^{-1}+ \mathnormal{Y}_{k}^{\mathsf{T}}\mathnormal{Y}_{k}\right)^{-1} \\
& \times \Phi^{\mathsf{T}}\left(k+N_{o}, k\right)+2 \beta_{2} \sum_{i=0}^{N_{o}-1} 2^{i} \\
& \times \Phi \left(k+N_{o}, k+N_{o}-i\right) \Phi^{\mathsf{T}}\left(k+N_{o}, k+N_{o}-i\right).
\end{aligned}
\end{equation}
Aggregating the bounds in Assumption~\ref{assum:bounds_matrix},~\eqref{eq:obsv_gramian},~\eqref{eq:Z_bound}, and~\eqref{eq:res22}, we have
\begin{equation}
\label{eq:int2prime}
    Y_{k}^{\mathsf{T}} Y_k \succeq \gamma_2 W_o (k+N_o,k),
\end{equation}
with $\gamma_2 = \frac{1}{\overline{\rho}} \prod_{j=k}^{k+N_o-1} \alpha_{2,j+1}$. Finally, we assimilate the inequalities in~\eqref{eq:int1prime} and~\eqref{eq:int2prime}, along with Assumptions~\ref{assum:agc_bounds} and~\ref{assum:obsv_bounds}, which gives us
\begin{equation}
    P_{k+N_o} \preceq \underbrace{\left( \frac{2^{N_o}}{\gamma_2 \varepsilon} + 2 \beta_2 \sum_{i=0}^{N_o-1} \left( 2 \overline{\alpha} \right)^i \right)}_{\overline{\pi}(N_o)} I.
\end{equation}
\qed

\emph{Proof of Theorem~\ref{thm:main_thm_iss}:} From the equations~\eqref{eq:v_app1} and~\eqref{eq:opt_sol}, we can obtain the dynamics of the estimation error $e[k]$, which admits the following form
\begin{equation}
    e[k+1] = (I - K_{k+1} C_{k+1}) (\tilde{A}_{\mathfrak{v}} e[k] - \tilde{G}_{\mathfrak{v}} r[k]) + K_{k+1} v[k+1].
\end{equation}
In order to prove exponential input-to-state stability of the estimation error, we consider the candidate Lyapunov function
\begin{equation}
    V_k = e[k]^{\mathsf{T}} P_k^{-1} e[k].
\end{equation}
Consider any time index $k$ that satisfies $k \geq \max \{ N_c,N_o \}$ and let $S_{k+1}$ be an arbitrary matrix. We have
\begin{align}
\label{eq:ISS1}
    &e[k+1]^{\mathsf{T}}\left(\mathnormal{P}_{k+1}+\mathnormal{S}_{k+1} \mathnormal{S}_{k+1}^{\mathsf{T}}\right)^{-1} e[k+1] \nonumber \\
    &=\left(\tilde{A}_{\mathfrak{v}} e[k]-\tilde{G}_{\mathfrak{v}} r[k]\right)^{\mathsf{T}} \mathnormal{M}_{k+1}^{-1}\left(\tilde{A}_{\mathfrak{v}} e[k]-\tilde{G}_{\mathfrak{v}} r[k]\right) \nonumber \\ 
    &+v[k+1]^{\mathsf{T}} \mathnormal{R}_{k+1}^{-1} v[k+1] \nonumber \\
    &-\left(C_{k+1}\left(\tilde{A}_{\mathfrak{v}} e[k]-\tilde{G}_{\mathfrak{v}} r[k]\right)-v[k+1]\right)^{\mathsf{T}} \nonumber \\
    & \times \left(C_{k+1} \mathnormal{M}_{k+1} C_{k+1}^{\mathsf{T}}+\mathnormal{R}_{k+1}\right)^{-1} \nonumber \\
    & \times \left(C_{k+1}\left(\tilde{A}_{\mathfrak{v}} e[k]-\tilde{G}_{\mathfrak{v}} r[k]\right)-v[k+1]\right) \nonumber \\
    &-\left(\mathnormal{M}_{k+1}^{-1}\left(\tilde{A}_{\mathfrak{v}} e[k]-\tilde{G}_{\mathfrak{v}} r[k]\right)+C_{k+1}^{\mathsf{T}} \mathnormal{R}_{k+1}^{-1} v[k+1]\right)^{\mathsf{T}} \mathnormal{S}_{k+1} \nonumber \\
    &\times \left(\mathnormal{I}+\mathnormal{S}_{k+1}^{\mathsf{T}}\left(\mathnormal{M}_{k+1}^{-1}+C_{k+1}^{\mathsf{T}} \mathnormal{R}_{k+1}^{-1} C_{k+1}\right) \mathnormal{S}_{k+1}\right)^{-1} \mathnormal{S}_{k+1}^{\mathsf{T}} \nonumber \\
    &\times\left(\mathnormal{M}_{k+1}^{-1}\left(\tilde{A}_{\mathfrak{v}} e[k]-\tilde{G}_{\mathfrak{v}} r[k]\right)+C_{k+1}^{\mathsf{T}} \mathnormal{R}_{k+1}^{-1} v[k+1]\right) \nonumber \\
    & \stackrel{(\bullet)}{\leq} \left(\tilde{A}_{\mathfrak{v}} e[k]-\tilde{G}_{\mathfrak{v}} r[k]\right)^{\mathsf{T}} \mathnormal{M}_{k+1}^{-1}\left(\tilde{A}_{\mathfrak{v}} e[k]-\tilde{G}_{\mathfrak{v}} r[k]\right) \nonumber \\
    &+2 v[k+1]^{\mathsf{T}} \mathnormal{R}_{k+1}^{-1} v[k+1] -\frac{1}{2}\left(\tilde{A}_{\mathfrak{v}} e[k]-\tilde{G}_{\mathfrak{v}} r[k]\right)^{\mathsf{T}} \nonumber \\
    & \times \mathnormal{M}_{k+1}^{-1} \mathnormal{S}_{k+1}\left(\mathnormal{I}+\mathnormal{S}_{k+1}^{\mathsf{T}}\left(\mathnormal{M}_{k+1}^{-1}+C_{k+1}^{\mathsf{T}} \mathnormal{R}_{k+1}^{-1} C_{k+1}\right) \mathnormal{S}_{k+1}\right)^{-1} \nonumber \\
    & \times \mathnormal{S}_{k+1}^{\mathsf{T}} \mathnormal{M}_{k+1}^{-1}\left(\tilde{A}_{\mathfrak{v}} e[k]-\tilde{G}_{\mathfrak{v}} r[k]\right) \nonumber \\
    & \stackrel{(\bullet \bullet)}{\leq} \left(1-\frac{\alpha_{3, k+1}}{2}\right)\left(\tilde{A}_{\mathfrak{v}} e[k]-\tilde{G}_{\mathfrak{v}} r[k]\right)^{\mathsf{T}} \mathnormal{M}_{k+1}^{-1} \nonumber \\
    & \times \left(\tilde{A}_{\mathfrak{v}} e[k]-\tilde{G}_{\mathfrak{v}} r[k]\right) +\frac{\alpha_{3, k+1}}{2}\left(\tilde{A}_{\mathfrak{v}} e[k]-\tilde{G}_{\mathfrak{v}} r[k]\right)^{\mathsf{T}} \nonumber \\
    &\times \left(\mathnormal{M}_{k+1}+\mathnormal{S}_{k+1} \mathnormal{S}_{k+1}^{\mathsf{T}}\right)^{-1}\left(\tilde{A}_{\mathfrak{v}} e[k]-\tilde{G}_{\mathfrak{v}} r[k]\right)\nonumber \\
    &+2 v[k+1]^{\mathsf{T}} \mathnormal{R}_{k+1}^{-1} v[k+1] \nonumber \\
    & \stackrel{(\square)}{\leq} \left(1-\frac{\alpha_{3, k+1}}{2}\right)\left(1+\varepsilon_{3}\right) e[k]^{\mathsf{T}} \mathnormal{P}_{k}^{-1} e[k]\nonumber\\
    &+\frac{\alpha_{3, k+1}}{2}\left(1+\varepsilon_{3}\right) e[k]^{\mathsf{T}} \nonumber \\
    & \times \left(\mathnormal{P}_{k}+\tilde{A}_{\mathfrak{v}}^{-1}\left(\tilde{G}_{\mathfrak{v}} \mathnormal{Q}_{k} \tilde{G}_{\mathfrak{v}}^{\mathsf{T}}+\mathnormal{S}_{k+1} \mathnormal{S}_{k+1}^{\mathsf{T}}\right) \tilde{A}_{\mathfrak{v}}^{-\mathsf{T}}\right)^{-1} e[k] \nonumber \\
    &+\left( \frac{1+\varepsilon_{3}}{\varepsilon_{3}} \right) r[k]^{\mathsf{T}} \mathnormal{Q}_{k}^{-1} r[k]+2 v[k+1]^{\mathsf{T}} \mathnormal{R}_{k+1}^{-1} v[k+1],
\end{align}
with $\alpha_{3,k+1} = \| I + \frac{1}{\underline{\pi}} S_{k+1}^{\mathsf{T}} S_{k+1} \|^{-1}$ and $\varepsilon_3 \in \mathbb{R}^{+}$. The inequality $(\bullet)$ in~\eqref{eq:ISS1} is a consequence of the Young-like inequality~\eqref{eq:younglike}, the inequality $(\bullet \bullet)$ in~\eqref{eq:ISS1} results from the Woodbury identity used in conjunction with Lemma~\ref{lemma:P1} and~\eqref{eq:mat_inv_P_inv}, whereas the inequality $(\square)$ is a result of~\eqref{eq:younglike} and~\eqref{eq:filter_updates}.

From Assumption~\ref{assum:bounds_matrix} and~\eqref{eq:ISS1}, we can write
\begin{equation}
\label{eq:ISS4}
\begin{aligned}
    &e[k+1]^{\mathsf{T}} P_{k+1}^{-1} e[k+1] \leq (1+\varepsilon_3) e[k]^{\mathsf{T}} P_{k}^{-1} e[k] \\
    &+ \left( \frac{1+\varepsilon_3}{\varepsilon_3 \underline{\vartheta}} \right) \| r[k] \|^2 + \left( \frac{2}{\underline{\rho}} \right) \| v[k+1] \|^2.
\end{aligned}
\end{equation}
Now, from the equality 
\begin{equation}
\label{eq:ISS3}
    S_k S_k^{\mathsf{T}} = \tilde{A}_{\mathfrak{v}}^{-1} ( \tilde{G}_{\mathfrak{v}} Q_k \tilde{G}_{\mathfrak{v}}^{\mathsf{T}} + S_{k+1} S_{k+1}^{\mathsf{T}} ) \tilde{A}_{\mathfrak{v}}^{-\mathsf{T}}    
\end{equation}
we have,
\begin{align}
\label{eq:ISS2}
    &e[k+1]^{\mathsf{T}}\left(\mathnormal{P}_{k+1}+\mathnormal{S}_{k+1} \mathnormal{S}_{k+1}^{\mathsf{T}}\right)^{-1} e[k+1] \nonumber \\
    &\leq \left( \frac{1+\varepsilon_3}{\varepsilon_3 \underline{\vartheta}} \right) \| r[k] \|^2 + \left(1-\frac{\alpha_{3, k+1}}{2}\right)\left(1+\varepsilon_{3}\right) e[k]^{\mathsf{T}} \mathnormal{P}_{k}^{-1} e[k] \nonumber \\
    &+ \left( \frac{2}{\underline{\rho}} \right) \| v[k+1] \|^2 \nonumber \\
    &+ \frac{\alpha_{3,k+1}}{2} (1 + \varepsilon_3) e[k]^{\mathsf{T}}\left(\mathnormal{P}_{k}+\mathnormal{S}_{k} \mathnormal{S}_{k}^{\mathsf{T}}\right)^{-1} e[k].
\end{align}
We let $S_{k+N_c} S_{k+N_c}^{\mathsf{T}} = 0$ and by repeated application of~\eqref{eq:ISS2}, we get
\begin{align}
\label{eq:ISS5}
    &e[k+N_c]^{\mathsf{T}} P_{k+N_c}^{-1} e[k+N_c] \nonumber \\
    &\leq (1 - \gamma_3)(1+\varepsilon_3)^{N_c} e[k]^{\mathsf{T}} P_k^{-1} e[k] \nonumber \\
    &+ \gamma_3 (1+\varepsilon_3)^{N_c} e[k]^{\mathsf{T}}\left(\mathnormal{P}_{k}+\mathnormal{S}_{k} \mathnormal{S}_{k}^{\mathsf{T}}\right)^{-1} e[k] \nonumber \\
    &+ \frac{(1+\varepsilon_3)^{N_c}}{\varepsilon_3 \underline{\vartheta}} \sum_{i=k}^{k+N_c-1} \| r[i] \|^2 \nonumber \\
    &+ \frac{2(1+\varepsilon_3)^{N_c-1}}{\underline{\rho}} \sum_{j=k}^{k+N_c-1} \| v[j+1] \|^2,
\end{align}
where $\gamma_3 = \frac{1}{2^{N_c}} \prod_{i=k}^{k+N_c-1} \alpha_{3,i+1}$. Using the recursive definition in~\eqref{eq:ISS3}, the bounds in~\eqref{eq:contr_gramian}, and Assumption~\ref{assum:bounds_matrix}, we get the bound
\begin{equation}
\label{eq:Sk_ineq}
    S_k S_k^{\mathsf{T}} \succeq \underline{\vartheta} \Phi^{-1} (k+N_c,k) W_c (k+N_c,k) \Phi^{-\mathsf{T}} (k+N_c,k).
\end{equation}
With the inequality~\eqref{eq:Sk_ineq}, we can then aggregate the bounds in Assumptions~\ref{assum:agc_bounds} and~\ref{assum:contr_bounds} and Lemma~\ref{lemma:P2} to then obtain
\begin{equation}
\label{eq:ISS6}
    S_k S_k^{\mathsf{T}} \succeq \frac{\underline{\vartheta} \delta}{\overline{\alpha}^{N_c} \overline{\pi}} P_k.
\end{equation}
Now, given the Lyapunov function $V_k = e[k]^{\mathsf{T}} P_k^{-1} e[k]$ and the bound in~\eqref{eq:ISS4}, we have
\begin{equation}
\label{eq:ISS7}
\begin{aligned}
        V_k &\leq (1+\varepsilon_3)^{k - k_0} V_{k_0} + \frac{(1+\varepsilon_3)^{k - k_0}}{\varepsilon_3 \underline{\vartheta}} \sum_{i=k_0}^{k-1} \| r[i] \|^2 \\
        &+ \frac{2(1+\varepsilon_3)^{k-k_0-1}}{\underline{\rho}} \sum_{j=k_0}^{k-1} \| v[j+1] \|^2
\end{aligned}
\end{equation}
for all $k \geq k_0 \geq \max \{ N_c,N_o \}$. Using Assumption~\ref{assum:bounds_matrix},~\eqref{eq:ISS5}, and~\eqref{eq:ISS6}, we have
\begin{equation}
\label{eq:ISS8}
\begin{aligned}
    V_{k+N_c} &\leq \eta_3 V_k + \frac{(1+\varepsilon_3)^{N_c}}{\varepsilon_3 \underline{\vartheta}} \sum_{i=k}^{k+N_c-1} \| r[i] \|^2 \\
    &+ \frac{2(1+\varepsilon_3)^{N_c-1}}{\underline{\rho}} \sum_{j=k}^{k+N_c-1} \| v[j+1] \|^2,
\end{aligned}
\end{equation}
for all $k \geq k_0 \geq \max \{ N_c,N_o \}$ and with
\begin{equation}
    \eta_3 = \left( 1 - \frac{\gamma_3 \underline{\vartheta} \delta}{\underline{\vartheta} \delta + \overline{\alpha}^{N_c} \overline{\pi}} \right) (1+\varepsilon_3)^{N_c}.
\end{equation}
If we assume, without loss of generality, that $\varepsilon_3$ is chosen such that $\eta_3 < 1$, then from~\eqref{eq:ISS7} and~\eqref{eq:ISS8}, we obtain
\begin{equation}
\label{eq:ISS9}
\begin{aligned}
    V_{k} &\leq \left(\frac{\left(1+\varepsilon_{3}\right)^{N_{c}}}{\eta_{3}}\right)^{\frac{N_{c}-1}{N_{c}}} \eta_{3}^{\frac{k-k_{0}}{N_{c}}} V_{k_{0}} \\ &+\frac{N_{c}\left(1+\varepsilon_{3}\right)^{N_{c}}}{\varepsilon_{3} \underline{\vartheta} \left(1-\eta_{3}\right)} \max _{k_{0} \leq i \leq k-1}\left\|r[i]\right\|^{2} \\
    &+\frac{2 N_{c}\left(1+\varepsilon_{3}\right)^{N_{c}-1}}{\underline{\rho} \left(1-\eta_{3}\right)} \max _{k_{0} \leq j \leq k-1}\left\| v[j+1] \right\|^{2}
\end{aligned}
\end{equation}
for all $k \geq k_0 \geq \max \{ N_c,N_o \}$.

On the other hand, from Lemmas~\ref{lemma:P1} and~\ref{lemma:P2}, we have the following bounds on the Lyapunov function
\begin{equation}
\label{eq:ISS10}
     \frac{1}{\overline{\pi}} \| e[k] \|^2 \leq V_k \leq \frac{1}{\underline{\pi}} \| e[k] \|^2,
\end{equation}
for all $k \geq \max \{ N_c,N_o \}$.

Thus, using~\eqref{eq:ISS9} and~\eqref{eq:ISS10}, the proof of the theorem follows with
\begin{subequations}
\begin{equation}
    \sigma = \sqrt{\frac{3 \overline{\pi}}{\underline{\pi}}} \left( \frac{(1+\varepsilon_3)^{N_c}}{\eta_3} \right)^{\frac{N_c-1}{2N_c}},
\end{equation}
\begin{equation}
    \tau = \eta_3^{\frac{1}{2N_c}},
\end{equation}
\begin{equation}
    \chi = \sqrt{\frac{3 \overline{\pi} N_c (1+\varepsilon_3)^{N_c}}{\varepsilon_3 \underline{\vartheta} (1-\eta_3)}},
\end{equation}
and
\begin{equation}
    \psi = \sqrt{\frac{6 \overline{\pi} N_c (1+\varepsilon_3)^{N_c-1}}{\underline{\rho} (1-\eta_3)}}.
\end{equation}
\end{subequations}
\qed

    




\bibliographystyle{IEEEtran}
\bibliography{IEEEabrv,mybibfile}

\end{document}